\documentclass[12pt]{amsart}
\usepackage{amscd,amssymb}
\usepackage[arrow,matrix]{xy}

\theoremstyle{plain}
\newtheorem{thm}[subsection]{Theorem}

\newtheorem{prop}[subsection]{Proposition}
\newtheorem{cor}[subsection]{Corollary}

\theoremstyle{definition}

\newtheorem{definition}[subsection]{Definition}
\newtheorem{ex}[subsection]{Example}

\numberwithin{equation}{section}
\newcommand{\A}{{\mathcal A}}

\newcommand{\CC}{{\mathcal C}}
\newcommand{\LL}{{\mathcal L}}

\newcommand{\V}{{\mathcal V}}

\newcommand{\al}{{\alpha}}

\newcommand{\Z}{\mathbb{Z}}

\newcommand{\C}{\mathbb{C}}

\newcommand{\PP}{\mathbb{P}}

\newcommand{\T}{\mathbb{T}}

\DeclareMathOperator{\Hom}{Hom}

\begin{document}

\title [  Admissible local systems for a class of line arrangements]
{ Admissible local systems for a class of line arrangements }

\author{SHAHEEN NAZIR, ZAHID RAZA}
 \address{School of Mathematical Sciences,
         Government College University,
         68-B New Muslim Town Lahore,
         PAKISTAN.}
\email {shaheen.nazeer@gmail.com}
 \address{Zahid Raza\\ School of Mathematical Sciences,
         Government College University,
         68-B New Muslim Town Lahore,
         PAKISTAN.}
\email {zahidsms@gmail.com }



\begin{abstract}
A rank one local system $\LL$ on a smooth complex algebraic variety $M$ is admissible roughly speaking if the dimension of the cohomology groups $H^m(M,\LL)$
can be computed directly from the cohomology algebra $H^*(M,\C)$.

We say that a line arrangement $\A$ is of type $\CC_k$ if $k \ge 0 $ is the minimal number of lines in $\A$
containing all the points of multiplicity at least 3.
We show that if $\A$ is a line arrangement in the classes $\CC_k$ for $k\leq 2$, then any rank one local system $\LL$ on the line arrangement complement $M$
is admissible. Partial results are obtained for the class $\CC_3$.

\end{abstract}

\maketitle

\section{Introduction } \label{s0}

When $M$ is a hyperplane arrangement complement in some projective space $\PP^n$, one  defines the notion of
an {\it admissible} local system $\LL$ on $M$ in terms of some conditions on the residues of an associated
logarithmic connection $\nabla(\al)$ on a good compactification of $M$, see for instance \cite{ESV},\cite{STV}, \cite{F}, \cite{LY}, \cite{DM}.
This notion plays a key role in the theory, since for such an admissible local system $\LL$ on $M$ one has
$$ \dim H^i(M,\LL)=\dim H^i(H^*(M,\C), \alpha \wedge)$$
for all $i$.

Let $\A$ be a line arrangement in the complex projective plane $\PP^2$ and denote by $M$ the corresponding
arrangement complement.
For the case of line arrangements, a good compactification as above is obtained just by blowing-up the points of multiplicity at least 3 in $\A$. This explains the simple version of the admissibility definition given below, see Definition \ref{d1}.

We say that a line arrangement $\A$ is of type $\CC_k$ if $k \ge 0$ is the minimal number of lines in $\A$
containing all the points of multiplicity at least 3. For instance $k=0$ corresponds to nodal arrangements, while $k=1$ corresponds to the case of a nodal {\it affine} arrangement, see \cite{CDP}.
Note that $k=k(\A)$ is combinatorially defined, i.e. depends only on the associated lattice $L(\A)$.

The main result of this note is the following one, which, as explained above, is new in the case $n=2$.

\begin{thm} \label{t0}

If $\A$ is a line arrangement in the classes $\CC_k$ for $k\leq 2$, then any rank one local system $\LL$ on $M$
is admissible.

\end{thm}

This result implies in particular that, for a combinatorially defined class of arrangements, the characteristic varieties are combinatorially determined by  $L(\A)$, see Theorem \ref{t3}.

In section 2 we explicite first the admissibility condition in the case of line arrangements.
Then we prove our Theorem and derive some consequenses for the characteristic varieties of
line arrangements in the classes $\CC_k$ for $k\leq 2$.

In the final section we find some sufficient conditions on a local system $\LL$ on the complement of
a line arrangement in the class $\CC_3$ which imply that $\LL$ is admissible, see Proposition \ref{prop7}.

The deleted $B_3$-arrangement discovered by A. Suciu shows that Theorem \ref{t0} does not hold for $k=3$,
so our result is the best possible. A discussion of this example from the point of view of our paper is given in Example \ref{ex0}. Further examples conclude the paper.

\section{Admissible rank one local systems} \label{s1}

Let $\A=\{L_0,L_1,...,L_n\}$ be a line arrangement in $\PP^2$ and set $M=\PP^2 \setminus (L_0 \cup...\cup L_n)$. Let $\T(M)=\Hom(\pi_1(M),\C^*)$ be the character group of
$M$. This is an algebraic torus
$\T(M) \simeq (\C^*)^{n}$. Consider the exponential mapping
\begin{equation}
\label{e1}
\exp :H^1(M,\C) \to H^1(M,\C^*)=\T(M)
\end{equation}
induced by the usual exponential function $ \C \to \C^*, ~~ t \mapsto \exp(2 \pi it)$, where $i=\sqrt {-1}$. Clearly $\exp(H^1(M,\C))=\T(M)$
and  $\exp(H^1(M,\Z))=\{1\}$.

A rank one local system $\LL\in \T(M)$ corresponds to the choice of some monodromy complex numbers
$\lambda _j \in \C^*$ for $0 \leq j \leq n$ such that $\lambda _0 ...\lambda _n=1$. And a cohomology class
$\alpha \in H^1(M,\C)$ is given by
\begin{equation}
\label{e2}
\alpha=\sum_{j=0,n}a_j\frac {df_j}{f_j}
\end{equation}
with the residues $a_j \in \C$ satisfying $\sum_{j=0,n}a_j=0$ and $f_j=0$ a linear equation for the line $L_j$. With this notation one has
$\exp (\alpha)=\LL$ if and only if $\lambda _j =\exp(2\pi i a_j)$ for any $j=0,...,n$.

\begin{definition} \label{d1}
A local system $\LL \in \T(M)$ as above is admissible if there is a cohomology class $\alpha \in H^1(M,\C)$ such that $\exp(\alpha)=\LL$, $a_j \notin \Z_{>0}$  and, for any point $p \in L_0 \cup...\cup L_n$ of multiplicity at least 3, one has
$$a(p)=\sum_ja_j \notin \Z_{>0}$$
where the sum is over all $j$'s such that $p \in L_j$. Set $b(p)=Re(a(p))$, where $Re$ denotes the real part of a complex number.

\end{definition}

\subsection{Proof of Theorem \ref{t0} } \label{s12}


We can assume that, for a line arrangement $\A$  of type $\CC_k$, the  lines in $\A$
containing all the points of multiplicity at least 3 are $L_0$,...,$L_{k-1}$.

 Given $\LL$, we can choose $Re(a_j) \in [0,1)$ for all
$j \geq 1$. Then the relation  $\sum_{j=0,n}a_j=0$ implies that
$$Re(a_0)=-\sum_{j=1,n}Re(a_j) \leq 0.$$
This completes the proof in the case when $\A$ is in the class $\CC_0$.

\medskip

Consider now the case when $\A$ is in the class $\CC_1$ and choose $p \in L_0$ as in Definition \ref{d1}.

If the set of lines passing through $p$ is $L_0, L_{j_1},...,L_{j_k}$ for some $k \geq 2$, then one has
$$b(p)=\sum_jRe(a_j)=-\sum_{k}Re(a_k) \leq 0$$
where the second sum is over all $k>0$ such that $p \notin L_k$. This shows that any local system $\LL \in \T(M)$
is admissible in this case.

\medskip

The case when $\A$ is in the class $\CC_2$ is much more subtle.  Let $m_1$ be the maximum of the real numbers
$b(p)$ where $p$ is a point of multiplicity at least 3 on the line $L_1$ but not on $L_0$ such that
$b(p) \in \Z_{>0}$. By convention we set $\max \emptyset =0$.

If $m_1=0$, then the initial choice is good, exactly as in the previous situation.

Assume now that  $m_1>0$ and let $p_1 \in L_1$ be a multiple point with $b(p_1)=m_1$.
In general this value $m_1$ can be attained at several points $p_1,~~p_1',~~p_1'',...$ on $L_1$.
Any such point will be called an {\it extremal point for $\LL$ on $L_1$}.

 Then we replace $a_1$ by $a_1'=a_1-m_1$ and hence $Re(a_1') \le 0$ if $Re(a_1')$ is an integer. We also replace $a_0$ by $a_0'=a_0+m_1$, such that
the total sum of the residues is still zero. Note that $Re(a_0') \leq 0$.

After these two changes, it is clear that $b(q) \notin \Z_{>0}$ for any point $q$ of multiplicity at least 3 on the line $L_1$ but not on $L_0$. If $O=L_0 \cap L_1$ is a point of multiplicity at least 3, then
$$b(O)=Re(a_0')+ Re(a_1')+\sum_{m}Re(a_m)= -\sum_{k}Re(a_k) \leq 0 $$
where the first sum is over $m=2,n$ such that $O \in L_m$ and the second sum is over all $k=2,n$ such that
$O \notin L_k$. Note that $O$ is never an extremal point.

If now $q \ne O$ is a point of multiplicity at least 3 on the line $L_0$, in the corresponding sum for
$b(q)$ we have the following terms

\medskip

(i) $-Re(a_j)$ for all $j>0$;

\medskip

(ii) $ +Re(a_{j_k})$ for all the lines different from $H_0$ passing through $q$;

\medskip

(iii)  $+Re(a_{m_k})$ for all the lines  passing through $p_1$.

\medskip

(in fact the terms of type (i) come from $Re(a_0)$, and the terms of type (iii) come from adding $m_1$ to
$a_0$ to get $a_0'$).

In the last two families there is at most one line in common, the line determined by $p_1$ and $q$.
It follows that the sum for $b(q)$ involves some negative terms from (i) and possibly one positive term
in the interval $[0,1)$. Therefore $b(q)<1$ and this concludes the proof of Theorem \ref{t0}.

\bigskip

The {\em characteristic varieties}\/ of $M$
are the jumping loci for the cohomology of $M$, with
coefficients in rank~$1$ local systems:
\begin{equation} \label{e3}
\V^i_k(M)=\{\rho \in \T(M) \mid \dim H^i(M, \LL_{\rho})\ge k\}.
\end{equation}
When $i=1$, we use the simpler notation $\V_k(M)=\V^1_k(M)$.
One has the following result, see \cite{Dadm}, Remark 2.9 (ii).

\begin{thm} \label{t3} If any local system in $\T(M)$ is admissible, then for any $k$ there are no translated components in the characteristic variety $\V_k(M)$.

In particular, for any irreducible component $W$ of some characteristic variety $\V_k(M)$, the dimension of
$H^1(M,\LL)$ is constant for $\LL \in W \setminus \{1\}$.

\end{thm}

Combining this result and our Theorem \ref{t0} we get the following.

\begin{cor} \label{c0}

If $\A$ is a line arrangement in the classes $\CC_k$ for $k\leq 2$, then for any $k$ there are no translated components in the characteristic variety $\V_k(M)$.

\end{cor}

\section{ Line arrangements in the class $\CC_3$  } \label{s3}

In this section we discuss only the situation when the 3 lines $L_0$, $L_1$ and $L_2$ containing all the
points of multiplicity at least 3 are concurrent. Let $O=L_0 \cap L_1 \cap L_2$.
Let $\LL \in \T(M)$ be a rank one local system and choose the residues $a_j$ as above.
Then we have a collection $P_1$ of extremal points for $\LL$ on $L_1$ and
a collection $P_2$ of extremal points for $\LL$ on $L_2$. By convention $P_j$ is empty if the corresponding maximum $m_j=0$ (see also the beginning of the proof below for the definition of $m_j$).
Since we start with the residues $a_j$'s such that $Re(a_j) \in [0,1)$ for all $j>0$, then we have as above $O \notin P_1 \cup P_2$.
For each point $q \in L_0$ of multiplicity at least 3, we denote by $\A_q$ the set of lines in $\A$
passing through $q$.

\begin{prop} \label{prop7}

With the above notation, assume that one of the following conditions hold

\medskip

\noindent (i) $P_1= \emptyset$;

\medskip

\noindent (ii) $P_2= \emptyset$;

\medskip

\noindent (iii)  for each point $q \in L_0$ of multiplicity at least 3, one has either

(a) $P_1 \setminus Y_q$ and $P_2 \setminus Y_q$ are both non-empty, where $Y_q$ is
the union of the lines in  $\A_q$, or

(b) one of the sets $P_1 \setminus Y_q$ and $P_2 \setminus Y_q$ is non-empty, say $p_1 \in P_1 \setminus Y_q$, and there is an extremal point $p_2 \in P_2$ such that the line determined by $p_1,p_2$ is not in $\A$.

\medskip

Then the local system $\LL$ is admissible.

\end{prop}

\proof

 Let $m_j$, for $j=1,2$, be the maximum of the real numbers
$b(p)$ where $p \ne O$ is a point of multiplicity at least 3 on the line $L_j$  such that
$b(p) \in \Z_{>0}$. The only essentially new case is when $m_1>0$ and $m_2>0$, i.e. the case (iii) above.
Then we set $a_j'=a_j-m_j$ for $j=1,2$ and $a_0'=a_0+m_1+m_2$.
Exactly as in the proof of Theorem \ref{t0}, this settles the case of multiple points on $L_1$ and $L_2$
distinct of $O$. At the point $O$ we get the following sum
$$b(O)=Re(a_0'+a_1'+a_2')+\sum_{m}Re(a_m)= -\sum_{k}Re(a_k) \leq 0 $$
where the first sum is over $m=3,n$ such that $O \in L_m$ and the second sum is over all $k=3,n$ such that
$O \notin L_k$.

If now $q \ne O$ is a point of multiplicity at least 3 on the line $L_0$, in the corresponding sum for
$b(q)$ we have the following terms

\medskip

(i) $-Re(a_j)$ for all $j>0$;

\medskip

(ii) $ +Re(a_{j_k})$ for all the lines different from $H_0$ passing through $q$;

\medskip

(iii)  $+Re(a_{m_k})$ for all the lines  passing through $p_1$, for a fixed choice $p_1 \in P_1$;

\medskip

(iv)  $+Re(a_{s_t})$ for all the lines  passing through $p_2$, for a fixed choice $p_2 \in P_2$.

\medskip

(in fact the terms of type (i) come from $Re(a_0)$, the terms of type (iii) come from adding $m_1$
and the terms of type (iv) come from adding $m_2$).

Our assumption (iii) says exactly that there is a choice of $p_1$
and $p_2$ such that in the last 3 families of positive terms there
is at most one line occuring twice (and none occuring three times).
It follows that the sum for $b(q)$ involves some negative terms from
(i) and possibly one positive term in the interval $[0,1)$.
Therefore $b(q)<1$ and this concludes the proof.

\endproof

\begin{ex} \label{ex0}
Here we analyse from the point of view of the above proof the non-admissible local systems in the case of
Suciu's arrangement, see \cite{S1}.

Choose: $L_0: x=0$, $L_1:x=y$, $L_2:y=0$ , $L_3: x=z$,
$L_4: y=z$, $L_5: x-y+z=0$, $L_6: x-y-z=0$, $L_7: z=0$ (this last one is the line at infinity).

There are $7$ points of multiplicity at least $3$, say, $O=L_0\cap
L_1\cap L_2=(0:0:1)$, $p_1=L_1\cap L_3\cap L_4=(1:1:1)$,
$p'_1=L_1\cap L_5\cap L_6\cap L_7=(1:1:0)$, $p_2=L_2\cap L_3\cap
L_6=(1:0:1)$, $p'_2=L_2\cap L_4\cap L_7=(1:0:0)$, $q_1=L_0\cap
L_4\cap L_5$,  $q_2=L_0\cap L_3\cap L_7=(0:1:0)$. It is clear that
$L_0,L_1,L_2$ contain all the points of multiplicity at lest 3, even
those at infinity.

Consider the local system $\LL_{\rho}$ which is used to define the
translated component
$$W=\rho\otimes \{(t,1,t^{-1},t^{-1},t,t^{-2},t^{2},1)\ |\ t\in \C^*\},$$
 where
$\rho=(1,-1,-1,-1,1,1,1,-1)$. The corresponding residues are
$$a_1=a_2=a_3=a_7=1/2, ~~ a_4=a_5=a_6=0 \text{ and } a_0=-2.$$
Using these residues we find out that $m_1=m_2=1$,
 $P_1=\{p_1,p'_1\}$
and $P_2=\{p_2,p'_2\}$.

For $q=q_2$, both sets $P_1\setminus Y_q$ and $P_2\setminus Y_q$ are empty.
 Thus the condition (iii) of Proposition \ref{prop7} is not fulfilled.
\end{ex}

\begin{ex} \label{ex0.5}

Here is an example of an arrangement in $\CC_3$ of the type discussed above for which
all local systems are admissible.
Let $\A$ be the arrangement in $\C^2$ consisting of the following 10 lines:

Choose $L_0: x=0$,  $L_1: y=1/2(x+3)$, $L_2: y=-1/2(x-3)$

$L_3: y=x+1$, $L_4:y=-(x-1)$,

$L_5: y=2(x+1)$, $L_6:y=-2(x-1)$,

$L_7: y=3/2(x+3)$, $L_8: y=-3/2(x+3)$,

$L_9: y=5/2(x-3)$ and $L_{10}: y=-5/2(x-3)$.

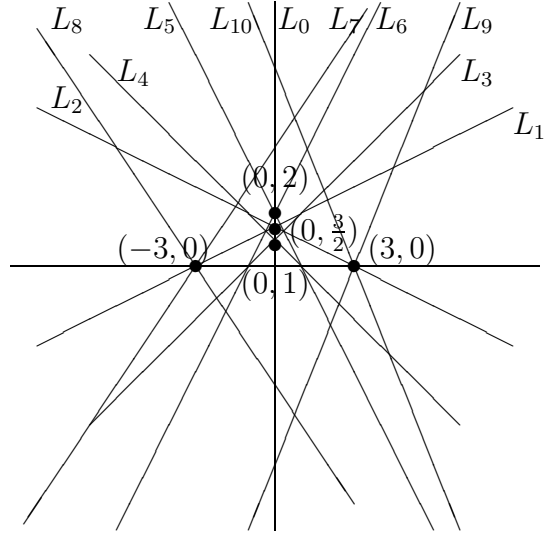
\vspace{20pt}
\begin{center}

 \begin{picture}(200,200)
 \put(100,120){\line(1,2){40}}
 \put(100,120){\line(-1,-2){60}}
 \put(100,145){\line(2,3){35}}
 \put(100,145){\line(-2,-3){95}}
 \put(100,55){\line(-2,3){90}}
  \put(100,55){\line(2,-3){30}}
 \put(100,25){\line(2,5){70}}
 \put(100,25){\line(-2,-5){10}}
 \put(100,175){\line(-2,5){10}}
 \put(100,175){\line(2,-5){70}}
 \put(100,120){\line(-1,2){40}}
  \put(100,120){\line(1,-2){60}}
 \put(100,110){\line(-1,1){70}}
 \put(100,110){\line(1,-1){70}}
 \put(100,110){\line(1,1){70}}
 \put(100,110){\line(-1,-1){70}}
 \put(100,115){\line(-2,1){90}}
 \put(100,115){\line(2,-1){90}}
 \put(100,115){\line(2,1){90}}
 \put(100,115){\line(-2,-1){90}}
 \put(100,0){\line(0,1){200}}
 \put(0,100){\line(1,0){200}}
 \put(97,117){$\bullet$}
 \put(97,105){$\bullet$}
 \put(97,111){$\bullet$}
 \put(67,97){$\bullet$}
 \put(127,97){$\bullet$}
 \put(101,190){$L_{0}$}
 \put(190,150){$L_{1}$}
 \put(15,160){$L_{2}$}
 \put(170,170){$L_{3}$}
 \put(40,170){$L_{4}$}
 \put(50,190){$L_{5}$}
 \put(138,190){$L_{6}$}
 \put(120,190){$L_{7}$}
 \put(15,190){$L_{8}$}
  \put(170,190){$L_{9}$}
   \put(75,190){$L_{10}$}

 \put(105,110){$(0,\frac{3}{2})$}
 \put(87,90){$(0,1)$}
 \put(87,130){$(0,2)$}
 \put(135,103){$(3,0)$}
 \put(40,103){$(-3,0)$}

 \end{picture}
 \end{center}

 \vspace{20pt}
\vspace{20pt}

There are 5 triple points, $3$ on the line $L_0$  namely $O=L_1\cap L_2\cap L_0=(0,\frac{3}{2})$,
 $q_1=L_0\cap L_3 \cap L_4 =(0,1)$ and $q_2=L_0\cap L_5 \cap L_6 =(0,2)$, the other 2 being
$p_1=L_1\cap L_7\cap L_8=(-3,0)$ and $p_2=L_2\cap L_9\cap L_{10}=(3,0)$.
Any local system is admissible by Proposition \ref{prop7}. Indeed, if we suppose that $P_1$ and $P_2$ are both non-empty, then
both sets $P_1\setminus Y_q=P_1$ and $P_2\setminus Y_q=P_2$ are again non- empty, for $q=q_1$ and $q=q_2$.

\end{ex}

\begin{ex} \label{ex1}
Here is an example of an arrangement in $\CC_3$ of the type discussed above for which
Proposition \ref{prop7} does not apply, but still most (perhaps all)
local systems are admissible.
Let $\A$ be the arrangement in $\C^2$ consisting of the following 7 lines:

$L_0: x=0$, $L_1: y=-2(x-1)$, $L_2: y=2(x+1)$, $L_4: y=x+1$, $L_5: y=-1/3(x+1)$,
$L_3: y=-(x-1)$, $L_6: y=1/3(x-1)$. The line at infinity is denoted by $L_7$.

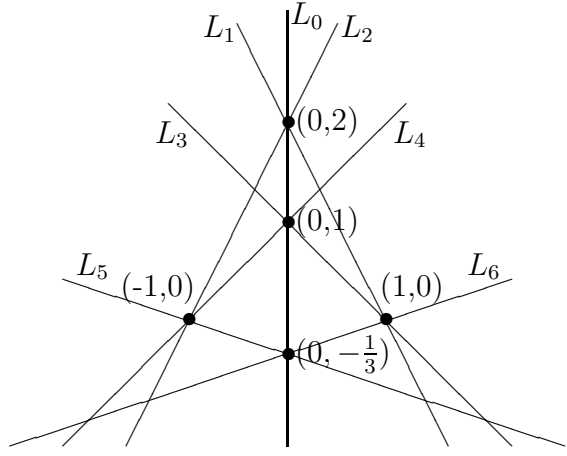
\vspace{20pt}

\vspace{20pt}
\begin{center}

 \begin{picture}(360,200)
 \put(104,0){\line(1,2){80}}
 \put(226,0){\line(-1,2){80}}
 \put(250,0){\line(-1,1){130}}
 \put(80,0){\line(1,1){130}}
 \put(60,0){\line(3,1){190}}
 \put(270,0){\line(-3,1){190}}
 \put(165,0){\line(0,1){165}}
 \put(162.5,120){$\bullet$}
 \put(162.5,82){$\bullet$}
 \put(162.5,32){$\bullet$}
 \put(125,45){$\bullet$}
 \put(199.5,45){$\bullet$}
 \put(166,160){$L_{0}$}
 \put(133,155){$L_{1}$}
 \put(185,155){$L_{2}$}
 \put(115,115){$L_{3}$}
 \put(205,115){$L_{4}$}
 \put(85,65){$L_{5}$}
 \put(233,65){$L_{6}$}
 \put(168,120){(0,2)}
 \put(168,82){(0,1)}
 \put(168,32){($0,-\frac{1}{3}$)}
 \put(102,57){(-1,0)}
 \put(200,57){(1,0)}
 \end{picture}
 \end{center}

 \vspace{20pt}

Here are $5$ points of multiplicity at least $3$ which are on the lines $L_0, L_1$ and $L_2$:

$O=L_0\cap \ L_1\cap L_2=(0,2)$, $p_1=L_1\cap \ L_3\cap L_6=(0,1)$,
$p_2=L_2\cap \ L_4\cap L_5=(-1,0)$,
$q_1=L_0\cap \ L_5\cap L_6=(0,-\frac{1}{3})$,
$q_2=L_0\cap \ L_3\cap L_4=(0,1)$.

Assume that we have a local system $\LL$ such that $P_1$ and $P_2$ are both non-empty.
Then one has
$$m_1=Re(a_2+a_4+a_5) \in \{1,2\}$$
and
$$m_2=Re(a_1+a_3+a_6) \in \{1,2\}.$$
In this case the new residue along $L_0$ is
$$a_0'=a_0+m_1+m_2= -\sum_{k=1,7}a_k+m_1+m_2=-a_7-i\cdot Im(\sum_{k=1,6}a_k)$$
where $Im$ denotes the imaginary part of a complex number.
When we compute the residue along the exeptional curve $E_1$ created by blowing-up the point $q_1$, we get
$$r_1=a_0'+a_5+a_6=a_5+a_6-a_7-i\cdot Im(\sum_{k=1,6}a_k).$$
Similarly we get
$$r_2=a_3+a_4-a_7-i\cdot Im(\sum_{k=1,6}a_k)$$
for the point $q_2$.
In this case Proposition \ref{prop7} does not apply, since
for $q=q_1,q_2$ both sets $P_1\setminus Y_q$ and $P_2\setminus Y_q$ are empty.
However, as soon as $r_j$ for $j=1,2$ are not  strictly positive integers, we know that the corresponding local system $\LL$ is admissible.

\end{ex}

\end{document}